\documentclass[a4paper, 12pt]{article}
 \usepackage{mathrsfs,amsfonts,amsmath,amssymb}
 \usepackage{color}

\makeatletter

\@addtoreset{equation}{section}
\makeatother
 \newtheorem{definition}{Definition}[section]
 \newtheorem{hypothesis}{Hypothesis}[section]
 \newtheorem{lemma}{Lemma}[section]
 \newtheorem{proposition}{Proposition}[section]
 \newtheorem{theorem}{Theorem}[section]
 \newtheorem{corollary}{Corollary}[section]
 \newtheorem{Remark}{Remark}[section]
 \newtheorem{example}{Example}[section]

 \def\beqlb{\begin{eqnarray}}\def\eeqlb{\end{eqnarray}}
 \def\beqnn{\begin{eqnarray*}}\def\eeqnn{\end{eqnarray*}}
 \def\beqn{\begin{equation}}
 \def\eeqn{\end{equation}}
\begin{document}

\title{\textbf{Large deviation principle of SDEs with non-Lipschitzian
coefficients under localized conditions}}
\author{Yunjiao Hu$^a$\quad and Guangqiang Lan$^b$\footnote{Corresponding author. Email: langq@mail.buct.edu.cn. Supported by China Scholarship Council,
National Natural Science Foundation of China (NSFC11026142) and Beijing Higher Education Young Elite Teacher Project (YETP0516).}
\\ \small $^{a,b}$ School of Science, Beijing University of Chemical Technology, Beijing 100029, China}

\date{}

\maketitle

\begin{abstract}
Localized sufficient conditions for the large deviation principle of the given stochastic differential
equations will be presented for stochastic differential equations with non-Lipschitzian and
time-inhomogeneous coefficients, which is weaker than those
relevant conditions existing in the literature. We consider at first the large
deviation principle when $\int_0^t\sup_{x\in\mathbb{R}^d}||\sigma(s,x)||\vee|b(s,x)|ds=:C_t<\infty$ for any fixed $t$,
then we generalize the conclusion to unbounded case by using bounded
approximation program.
 \end{abstract}



\section{ Introduction and Main Results}

Let $(\Omega,\mathscr{F},P)$ be a probability space, endowed with a
complete filtration $(\mathscr{F}_t)_{t\geq 0}$. Consider
the following stochastic differential equations (SDEs for short):
\beqn\label{SDE} dX^\varepsilon_t=\varepsilon^{\frac{1}{2}}\sigma(t,X^\varepsilon_t)dB_t
+b(t,X^\varepsilon_t)dt,\ X^\varepsilon_0=x_0\in \mathbb{R}^d
\eeqn
where $\varepsilon$ is an arbitrary positive number, $B_t$ is an
$m$-dimensional standard $\mathscr{F}_t$-Brownian motion, and
$\sigma$ and $b$ are $\mathscr{F}_t$-adapted functions from
$\mathbb{R}\times \mathbb{R}^d$ to
$\mathbb{R}^d\otimes\mathbb{R}^m$ and $\mathbb{R}^d$, respectively.
Without loss of generality, we restrict ourselves that $t$ is on the
interval [0,1]. Let $\mu_\varepsilon$ be the law of the solution of
(\ref{SDE}) on the space
$C_{x_0}([0,1],R^d)$ of continuous functions starting from $x_0\in
R^d.$ When the coefficients are time-homogeneous, Fang and Zhang
\cite{FZ} got the large deviation principle under certain global
non-Lipschitzian conditions. In \cite{Lan}, the second named author generalized Fang
and Zhang's result to a more general time-inhomogeneous case. However,
both of \cite{FZ,Lan} were concerned with global conditions. The aim of this
paper is to get the large deviation principle under localized conditions which are
even weaker than local Lipschitzian conditions (see the following (\ref{lianxu1})).

Much of the earlier work on the large deviation principle (see e.g. Freidlin
and Wentzell \cite{VF}, Donsker and Varadhan \cite{DV,DV1}) was based
on change of measure techniques, where a new measure is identified
under which the events of interest have high probability, and then
the probability of that event under the original probability measure
is calculated using the Radon-Nikodym derivative. An approach analogous
to the Prohorov compactness approach to weak
convergence has been developed by Pukhalskii \cite{Puk}, O'Brien and
Vervaat \cite{OV}, de Acosta \cite{Aco}. In \cite{CF}, Chiarini and Fischer got the sufficient
conditions of the so-called Laplace principle for stochastic differential
equations when the coefficients also depend on $\varepsilon$ and the past
of solution trajectory. The proof is based on the weak convergence approach
introduced by Dupuis and Ellis \cite{DE}. For more results about large
deviation principle, one can see \cite{Aco1,BC,BDM,Varadhan1},
and references therein.

In order that the integrals in the definition of the solutions of the equation (\ref{SDE}) are well-defined, we make
the following assumption which is enforced throughout the paper
\begin{equation}\label{keji}\int_0^T\sup_{|x|\le R}(|b(s,x)|+||\sigma(s,x)||^2)ds<\infty,\quad \forall T,R>0.
\end{equation}

We also assume throughout the paper that the coefficients $\sigma$ and $b$ satisfy
\beqn\label{lianxu} ||\sigma(t,x)-\sigma(t,y)||+|b(t,x)-b(t,y)|\le G(t)H(|x-y|),\eeqn
where $H:\mathbb{R}_+\rightarrow\mathbb{R}_+$ is an increasing continuous function with
$H(0)=0,$ and $0\le G\in L^2([0,t])$ for any fixed $t>0.$ $|\cdot|$ denotes the
Euclidean distance and $\langle\cdot,\cdot\rangle $
inner product in $\mathbb{R}^d$, $||\sigma||^2=\sum_{i=1}^d\sum_{j=1}^m\sigma_{ij}^2$.

We consider the sufficient conditions of large deviation principle of stochastic differential equation
(\ref{SDE}). Assume that for any fixed $R>0$, $\eta_R:[0,1)\rightarrow\mathbb{R}_+$ is a differentiable
function satisfying $\int_{0+}\frac{dx}{\eta_R(x)}=\infty$ with $\eta_R(0)=0, \eta_R'\ge0$,
and $f$ is nonnegative function such that $\int_0^t f^2(s)ds<\infty, \forall t>0.$

Then we have the following result

\begin{theorem}\label{Euler}
Let $R>0$ be fixed arbitrarily. Assume that 
\beqn\label{keji1}\int_0^t\sup_{x\in\mathbb{R}^d}\big(||\sigma(s,x)||
\vee|b(s,x)|\big)ds<\infty\eeqn
for any fixed $t>0$. If for any $|x|\vee|y|\le R,$ the
following condition
\beqn\label{lianxu1}\Big(||\sigma(t,x)-\sigma(t,y)||^2+2\langle x-y,
b(t,x)-b(t,y)\rangle\Big)\vee\Big|(\sigma(t,x)-\sigma(t,y))^T(x-y)\Big|
\leq f(t)\eta_R(|x-y|^2)
 \eeqn
holds with $|x-y|<c_0(<1), t\in [0,1],$ (here and from now on $A^T$ denotes the 
transpose of a matrix $A$), then the distribution family
$\{\mu_\varepsilon,\varepsilon>0\}$ satisfies a
large deviation principle with the following good rate function
$$I(h)=\inf\{\frac{1}{2}\int_0^1|l'(t)|^2dt:F(l)=h,\ l\in C_0([0,1],R^m)\},
\quad h\in C_{x_0}([0,1],R^d),$$
where $l'$ denotes the gradient of $l$, and $F(l)$ satisfies the auxiliary ordinary differential equation
\beqn\label{fuzhu} F(l)(t)=x_0+\int_0^tb(s,F(l)(s))ds+\int_0^t\sigma(s,F(l)(s))l'(s)ds \eeqn
for $l\in C_0([0,1],R^m),t>0.$ \end{theorem}

To prove this result, we need the following lemmas.

Suppose $X_n^\varepsilon(t)$ is Euler approximation of
$X^\varepsilon(t)$ defined as
\beqn\label{EMX}
X_n^{\varepsilon}(t):=x_0+\int_0^tb(s,X_n^{\varepsilon}(\frac{[ns]}{n}))ds
+\varepsilon^{\frac{1}{2}}\int_0^t\sigma(s,X_n^{\varepsilon}(\frac{[ns]}{n}))dB_s
\eeqn
where the $[ns]$ denotes the integral part of $ns$.

\begin{lemma}\label{l1} Assume that $\sigma$ and $b$ are bounded and that
(\ref{lianxu1}) holds. Then for any $\delta_0>0$, we have
$$
\lim_{n\rightarrow\infty}\limsup_{\varepsilon\rightarrow0}\
\varepsilon \log P(\sup_{0\le
t\le1}|X^\varepsilon(t)-X_n^\varepsilon(t)|\ge\delta_0)=-\infty .
$$
\end{lemma}

Let
$l\in C_0([0,1],R^m)$, define
\[e(l):=\left\{\begin{array}{ll} \int_0^1|l'(t)|^2dt,\quad
\textrm{if}\ l\ \textrm{is absolutely continuous},\\
+\infty,\qquad\qquad \textrm{otherwise}.
\end{array}
\right.\]
and $F_n$ be the Euler approximation of $F$, namely,
\beqn\label{EMF}
F_n(l)(t):=x_0+\int_0^tb(s,F_n(l)(\frac{[ns]}{n}))ds+\int_0^t\sigma(s,F_n(l)(\frac{[ns]}{n}))l'(s)ds.
\eeqn

\begin{lemma}\label{l2} Assume that $\sigma$ and $b$ are bounded and satisfy
\beqn\label{lianxu2}\langle x-y,b(t,x)-b(t,y)\rangle\vee\Big|(\sigma(t,x)-\sigma(t,y))^T(x-y)\Big|
\leq f(t)\eta_R(|x-y|^2)
\eeqn
where $\eta_R$ and $f$ are same as in condition (\ref{lianxu1}).
Then for any $\alpha>0$,
$$
\lim_{n\rightarrow\infty}\sup_{\{l:e(l)\le\alpha\}} \sup_{0\le
t\le1}|F_n(l)(t)-F(l)(t)|=0.
$$
\end{lemma}

Notice that condition (\ref{lianxu2}) is weaker than (\ref{lianxu1}).

From now on consider the general case that $\sigma$ and $b$ only satisfy the integrable condition (\ref{keji}).
Let $\gamma:\mathbb{R}_+\rightarrow\mathbb{R}_+$ is a differentiable function such that 
$\gamma'(x)\ge0,\ \lim_{s\rightarrow\infty}\gamma(s)=\infty$, and $\int_K^{\infty}\frac{ds}{\gamma(s)+1}=\infty$
holds for some $K>0$, $g$ is a nonnegative function such that $\int_0^t g^2(s)ds<\infty.$

\begin{lemma}\label{l3} If there exists $K>0$ such that
\beqn\label{zengzhang}\Big(||\sigma(t,x)||^2+2\langle x,b(t,x)\rangle\Big)\vee|\sigma^T(t,x)x| \leq
g(t)\big(\gamma(|x|^2)+1\big)
\eeqn
holds for any $|x|\ge K,\ t\in [0,1]$, then
$$
\lim_{R\rightarrow\infty}\limsup_{\varepsilon\rightarrow0}\
\varepsilon \log P(\sup_{0\le t\le1}|X^\varepsilon(t)|\ge R)=-\infty.
$$
\end{lemma}

\begin{lemma}\label{l4}
Let $\alpha>0.$ Assume that (\ref{zengzhang}) holds. Then
$$\sup_{\{l:e(l)\le\alpha\}} \sup_{0\le t\le1}|F(l)(t)|<\infty.$$
\end{lemma}

Note that, by the results we have proved in \cite{Lan1,Lan2}, if $\sigma$ and $b$
satisfy conditions (\ref{lianxu1}) and (\ref{zengzhang}), then for any fix $\varepsilon>0,$
stochastic differential equation (\ref{SDE}) has a unique global strong solution.

\begin{lemma}\label{l5} $I(\cdot)$ defined as in Theorem
\ref{Euler} is a good rate function on $C_{x_0}([0,1],R^d)$, that is, for any
$\alpha> 0,$ the sublevel set $\{h\in C_{x_0}([0,1],R^d):I(h)\le\alpha\}$ is compact.
\end{lemma}

\begin{theorem}\label{wujie} Suppose $\sigma$ and $b$ satisfy conditions (\ref{lianxu1}) 
and (\ref{zengzhang}). Then the distribution family $\{\mu_\varepsilon,\varepsilon>0\}$ 
of solutions of SDEs (\ref{SDE}) satisfy a large deviation principle with the same good 
rate function as in Theorem \ref{Euler}.\end{theorem}

\begin{Remark} In \cite{CF}, the authors got the large deviation principle
for the It\^o SDE under assumptions $H1$-$H6$, but the uniformly continuous condition
$H1$ does not hold in our case since $G(t)$ is only locally integrable but not bounded
 in $t$, and tightness condition $H6$ does not hold either since they used the sublinear 
 growth condition of $\sigma$ and $b$ to prove the tightness condition (see (\ref{lianxu1}) 
 and (\ref{zengzhang})). So their Theorem 2.1 can not be used under our condition. Moreover, 
 their assumptions are not easy to check. On the other hand, we can take 
 $\eta_R(x)=Rx\log\frac{1}{x},\ x\le\frac{1}{e}$ and $\gamma(x)=x\log x,\ x\ge K(>1)$ 
 for some $K$ large enough, $\sigma$ and $b$ may not satisfy the Lipschitz condition, 
 so the method of moments used in the literature  (such as \cite{DZ,DS,OP}) does not work 
 here because of the non-Lipschitz feature of coefficients.\end{Remark}

We now give an example to show that our conditions are really weaker than
those relevant conditions existing in the literature.

\begin{example}
We just consider the time-homogeneous case for simplicity.
Suppose $d=2,$ $m=1$. For any $r>0,$ define
$$\sigma(x)=|x|^r(-x_2,x_1)^T, b(x)=-|x|^{2r}x^T.$$

Since the local Lipschitzian condition holds for both
$\sigma$ and $b,$ the condition (\ref{lianxu1}) holds naturally. It's obvious that
there exists a unique strong solution for the giving stochastic
differential equation. On the other hand,
\[\Big(|\sigma(x)|^2+2\langle x,b(x)\rangle\Big)\vee |\sigma^T(x)x|
=(|x|^{2r+2}-2|x|^{2r+2})\vee 0=0\le K(|x|^2+1).\]

So by Theorem \ref{wujie}, we know that the large deviation principle holds in this case.
But there is NO constant $C>0$ such that
$$|\sigma(x)|^2=|x|^{2r+2}\le C(|x|^2\log |x|+1)$$
holds for $|x|$ large enough. So we have given a sufficient condition for large
deviation principle which is weaker than that of \cite{FZ}.
\end{example}

The rest of the paper is organized as follows. For the case that $\sigma(t,x)$ and $b(t,x)$ satisfy 
condition (\ref{keji1}), we will first show Lemma \ref{l1} in Section 2. Then we prove that the Euler 
approximation (\ref{EMF}) converges uniformly to the solution of auxiliary equation (\ref{fuzhu}) in 
Section 3. In Section 4, we will drop the assumption (\ref{keji1}) and get the large deviation principle 
by bounded approximation in general case.

\section{Proof of Lemma \ref{l1}}\vskip0.2in

Let $X_n^\varepsilon(t)$ be the Euler approximation of $X^\varepsilon(t)$ defined as (\ref{EMX}). Denote
$$Y_n^\varepsilon(t):=X_n^\varepsilon(t)-X^\varepsilon(t),\quad\xi_n^\varepsilon(t):=|Y_n^\varepsilon(t)|^2.$$

Define the following test function
$$
\varphi_{\rho,\lambda}(x):=\exp{\big(\lambda\int_0^x\frac{ds}{\eta_R(s)+\rho}\big)}.
$$

Then
$$
\varphi'_{\rho,\lambda}(x)=\frac{\lambda\varphi_{\rho,\lambda}(x)}{\eta_R(x)+\rho},$$

and
$$\varphi''_{\rho,\lambda} (x)=\frac{\lambda^2\varphi_{\rho,\lambda}
(x)-\lambda\varphi_{\rho,\lambda}(x)\eta_R'(x)}{(\eta_R(x)+\rho)^2}\le
\frac{\lambda^2\varphi_{\rho,\lambda} (x)}{(\eta_R(x)+\rho)^2}.
$$

We have used the condition $\eta'_R\ge0$ in the above inequality. Introduce stopping times
$$
\tau_n^\varepsilon:=\inf\{t>0,|X_n^\varepsilon(t)-X_n^\varepsilon(\frac{[nt]}{n})|\geq
\delta\}$$ 
where $\delta>0$ is an arbitrarily small number,
$$\tau_R:=\inf\{t>0,|X_n^\varepsilon(t)|\vee|X^\varepsilon(t)|\ge R\}$$

and
$$
T_n^\varepsilon:=\inf\{t>0,|\xi_n^\varepsilon(t)|\geq\delta_0^2\}.$$

Since $\sigma$ and $b$ satisfy integrable condition (\ref{keji1}), it's obvious that
$\lim_{R\rightarrow\infty}\tau_R=\infty.$

Without loss of generality, suppose $0<\delta_0\le c_0,$ then it
follows by It\^o's formula that
$$
\aligned\varphi_{\rho,\lambda}( \xi_n^\varepsilon(t)
)&=1+2\varepsilon^{\frac{1}{2}}\int_0^t\varphi'_{\rho,\lambda}(\xi_n^\varepsilon(s))
\langle Y_n^\varepsilon(s),(\sigma(s,X^\varepsilon(s))-\sigma(s,X_n^\varepsilon(\frac{[ns]}{n})))dB_s\rangle\\
&\quad+2\int_0^t\varphi'_{\rho,\lambda}( \xi_n^\varepsilon(s) )
\langle Y_n^\varepsilon(s),b(s,X^\varepsilon(s))-b(s,X_n^\varepsilon(\frac{[ns]}{n}))\rangle ds\\
&\quad+\varepsilon\int_0^t \varphi'_{\rho,\lambda}(
\xi_n^\varepsilon(s) )
||\sigma(s,X^\varepsilon(s))-\sigma(s,X_n^\varepsilon(\frac{[ns]}{n}))||^2 ds\\
&\quad +2\varepsilon\int_0^t\varphi''_{\rho,\lambda}(
\xi_n^\varepsilon(s)
)|(\sigma(s,X^\varepsilon(s))-\sigma(s,X_n^\varepsilon(\frac{[ns]}{n})))^TY_n^\varepsilon(s)|^2ds
.\endaligned $$

By definition of the stopping time, when $s\le\tau_n^\varepsilon\wedge T_n^\varepsilon\wedge\tau_R$,
it follows that $ |Y_n^\varepsilon(s)|\le\delta_0$ and $
|X_n^\varepsilon(s)-X_n^\varepsilon(\frac{[ns]}{n})|\leq \delta.$ So
\beqn\aligned&\langle
Y_n^\varepsilon(s),b(s,X^\varepsilon(s))-b(s,X_n^\varepsilon(\frac{[ns]}{n}))\rangle
\\&=\langle Y_n^\varepsilon(s),b(s,X^\varepsilon(s))-b(s,X_n^\varepsilon(s))\rangle+\langle Y_n^\varepsilon(s),b(s,X_n^\varepsilon(s))-b(s,X_n^\varepsilon(\frac{[ns]}{n}))\rangle\\&
\le \langle Y_n^\varepsilon(s),b(s,X^\varepsilon(s))-b(s,X_n^\varepsilon(s))\rangle+\delta_0G(s)H(\delta).\endaligned\eeqn

Similarly,
$$\aligned&||\sigma(s,X^\varepsilon(s))-\sigma(s,X_n^\varepsilon(\frac{[ns]}{n}))||^2
\\&\le2||\sigma(s,X^\varepsilon(s))-\sigma(s,X_n^\varepsilon(s))||^2+
2||\sigma(s,X_n^\varepsilon(s))-\sigma(s,X_n^\varepsilon(\frac{[ns]}{n}))||^2
\\&\le 2||\sigma(s,X^\varepsilon(s))-\sigma(s,X_n^\varepsilon(s))||^2+2G^2(s)H^2(\delta).\endaligned
$$

When $\varepsilon\le\frac{1}{2},$ by (\ref{lianxu1}) we arrive at
\[\aligned &2\langle
Y_n^\varepsilon(s),b(s,X^\varepsilon(s))-b(s,X_n^\varepsilon(\frac{[ns]}{n}))\rangle+
\varepsilon||\sigma(s,X^\varepsilon(s))-\sigma(s,X_n^\varepsilon(\frac{[ns]}{n}))||^2\\&
\le 2\langle Y_n^\varepsilon(s),b(s,X^\varepsilon(s))-b(s,X_n^\varepsilon(s))\rangle+2\varepsilon||\sigma(s,X^\varepsilon(s))-\sigma(s,X_n^\varepsilon(s))||^2\\&\quad+ 2G(s)H(\delta)\delta_0+2\varepsilon G^2(s)H^2(\delta)\\&\le f(s)\eta_R(\xi_n^\varepsilon(s))+2G(s)H(\delta)\delta_0+2\varepsilon G^2(s)H^2(\delta)\\&
\le \big(f(s)+2(G(s)+1)^2\big)[\eta_R(\xi_n^\varepsilon(s))+H(\delta)(\delta_0+\varepsilon H(\delta))].
\endaligned\]

On the other hand,
$$\aligned\Big|\big[\sigma(s,X^\varepsilon(s))-\sigma(s,X_n^\varepsilon(\frac{[ns]}{n}))\big]^TY_n^\varepsilon(s)\Big|^2
&\le 2\Big|(\sigma(s,X^\varepsilon(s))-\sigma(s,X_n^\varepsilon(s)))^TY_n^\varepsilon(s)\Big|^2\\&
\quad+2\Big|(\sigma(s,X_n^\varepsilon(s))-\sigma(s,X_n^\varepsilon(\frac{[ns]}{n})))^TY_n^\varepsilon(s)\Big|^2\\&
\le 2f^2(s)\eta_R^2(\xi_n^\varepsilon(s))+2\delta_0^2G^2(s)H^2(\delta)\\&
\le 2(f^2(s)+G^2(s))[\eta_R(\xi_n^\varepsilon(s))+\delta_0H(\delta)]^2.\endaligned$$

Take $$\rho=H(\delta)(\delta_0+\varepsilon H(\delta)).$$

Then $\rho\rightarrow 0$ as $\delta\rightarrow 0.$ By the definition of $\varphi_{\rho,\lambda}$, we have
$$\aligned \mathbb{E}\varphi_{\rho,\lambda}(\xi_n^\varepsilon
(t\wedge\tau_n^\varepsilon\wedge T_n^\varepsilon\wedge\tau_R))&\le
1+\lambda\mathbb{E}\int_0^{t\wedge\tau_n^\varepsilon\wedge T_n^\varepsilon\wedge\tau_R}(f(s)+2(G(s)+1)^2)
\varphi_{\rho,\lambda}(\xi_n^\varepsilon(s))ds\\&
\quad+4\varepsilon\lambda^2\mathbb{E}\int_0^{t\wedge\tau_n^\varepsilon\wedge T_n^\varepsilon\wedge\tau_R}(f^2(s)+G^2(s))
\varphi_{\rho,\lambda}(\xi_n^\varepsilon(s))ds. \endaligned$$

By Gronwall's lemma, it follows that (letting $R\rightarrow\infty$ and taking $t=1$),
$$ \mathbb{E} \varphi_{\rho,\lambda}(\xi_n^\varepsilon(1\wedge\tau_n^\varepsilon\wedge
T_n^\varepsilon))\le\exp\{\lambda\int_0^1(f+2(G+1)^2)(s)ds
+4\varepsilon\lambda^2\int_0^1(f^2+G^2)(s)ds\}.
$$

On the other hand,
$$ \aligned\mathbb{E} \varphi_{\rho,\lambda}(\xi_n^\varepsilon(1\wedge
\tau_n^\varepsilon\wedge T_n^\varepsilon))&\ge
 \mathbb{E}(\varphi_{\rho,\lambda}(\xi_n^\varepsilon(T_n^\varepsilon))I_{\{\tau_n^\varepsilon\ge1,T_n^\varepsilon\le 1\}}).
\\&=\varphi_{\rho,\lambda}(\delta_0^2)P(\tau_n^\varepsilon\ge1,T_n^\varepsilon\le
1).\endaligned
$$

Taking $\lambda=\frac{1}{\varepsilon} $. It follows that
$$
\varepsilon\log P(\tau_n^\varepsilon\ge1,T_n^\varepsilon\le 1)\le
C-\int_0^{\delta_0^2}\frac{ds}{\eta(s)+\rho},
$$
where $C:=\int_0^1(f(s)+4f^2(s))ds+\int_0^1(6G^2(s)+4G(s)+2)ds$. Therefore
$$
\limsup_{\varepsilon\rightarrow0}\ \varepsilon\log
P(\tau_n^\varepsilon\ge1,T_n^\varepsilon\le 1)\le
C-\int_0^{\delta_0^2}\frac{ds}{\eta(s)+\rho}.
$$

Now
$$
P(\sup_{0\le t\le1}|X_t^\varepsilon-X_n^\varepsilon(t)|\ge\delta_0)=P(T_n^\varepsilon\le
1)\le P(\tau_n^\varepsilon\ge1,T_n^\varepsilon\le 1)+P(\tau_n^\varepsilon\le1).
$$

Since $\sigma$ and $b$ satisfy integrable condition (\ref{keji1}), then
\[\aligned P(\tau_n^\varepsilon\le1)&\le\sum_{k=1}^nP(\sup_{\frac{k-1}{n}\le t<\frac{k}{n}}|X_n^\varepsilon(t)-X_n^\varepsilon(\frac{k-1}{n})|\ge\delta)\\&
\le\sum_{k=1}^nP(\sup_{\frac{k-1}{n}\le t<\frac{k}{n}}\varepsilon^\frac{1}{2}|
\int_{\frac{k-1}{n}}^t\sigma(s,X_n^\varepsilon
(\frac{k-1}{n}))dB_s|\ge\delta-B_{n,k}),\endaligned\]

\noindent where $B_{n,k}:=\int_{\frac{k-1}{n}}^\frac{k}{n}\sup_{x\in\mathbb{R^d}}|b(t,x)|dt, B:=\sum_{k=1}^nB_{n,k}.$
Similarly, define $A_{n,k}$ and $A$ with $|b|$ replaced by $||\sigma||^2.$ By (\ref{keji1}) again, it's obvious that
\beqn\label{jixian}\lim_{n\rightarrow\infty}\sup_{k\le n}(A_{n,k}\vee B_{n,k})=0.\eeqn 

Since
$$M_t:=\varepsilon^\frac{1}{2}\int_{\frac{k-1}{n}}^{t+
\frac{k-1}{n}}\sigma(s,X_n^\varepsilon(\frac{k-1}{n}))dB_s$$

\noindent is a martingale with respect to $\tilde{\mathscr{F}}_t:=\mathscr{F}_{t+\frac{k-1}{n}}$
for $0\le t<\frac{1}{n},$ then by exponential martingale inequality,
for any $d$ dimensional vector $|\theta|=1$, we have
\[\aligned P(\sup_{0\le t<\frac{1}{n}}\langle\theta,M_t\rangle\ge\delta-B_{n,k})&
\le P(\sup_{0\le t<\frac{1}{n}}\alpha\langle\theta,M_t\rangle-\frac{\alpha^2}{2}\big\langle\langle\theta,M\rangle\big\rangle_t
\ge\alpha(\delta-B_{n,k})-\frac{\varepsilon\alpha^2A_{n,k}}{2})\\&\le\exp(-\alpha(\delta
-B_{n,k})+\frac{\varepsilon\alpha^2A_{n,k}}{2}),\endaligned\]

\noindent where $\big\langle\langle\theta,M\rangle\big\rangle_t$ denotes the quadratic variation process
of $\langle\theta,M_t\rangle$. Taking $\alpha=\frac{\delta-B_{n,k}}{\varepsilon A_{n,k}}$,
\[\aligned P(\sup_{0\le t<\frac{1}{n}}
\langle\theta,M_t\rangle\ge\delta-B_{n,k})\le\exp(-\frac{\delta-B_{n,k}}{2\varepsilon A_{n,k}})\le\exp(-\frac{\delta-\sup_{k\le n}B_{n,k}}{2\varepsilon \sup_{k\le n}A_{n,k}}).\endaligned\]

We have used the fact (\ref{jixian}) here. Then by Stroock \cite{Stroock}, we have
\[\aligned P(\sup_{0\le t<\frac{1}{n}}
\varepsilon^\frac{1}{2}|\int_{\frac{k-1}{n}}^t\sigma(s,X_n^\varepsilon(\frac{k-1}{n}))dB_s|\ge\delta-B_{n,k})\le2d
\exp(-\frac{\delta-\sqrt{d}\sup_{k\le n}B_{n,k}}{2\varepsilon d\sup_{k\le n}A_{n,k}}).\endaligned\]

Thus,
\[\aligned P(\tau_n^\varepsilon\le1)\le2nd
\exp(-\frac{\delta-\sqrt{d}\sup_{k\le n}B_{n,k}}{2\varepsilon d\sup_{k\le n}A_{n,k}}).\endaligned\]

For sufficiently large $n$, it follow that
$$\aligned\varepsilon\log
P(\tau_n^\varepsilon\le1)\le\varepsilon\log(2nd)-\frac{\delta-\sqrt{d}\sup_{k\le n}B_{n,k}}{2d\sup_{k\le n}A_{n,k}}\le-\frac{\delta}{4d\sup_{k\le n}A_{n,k}}.\endaligned $$

Since $n$ is independent of $\varepsilon,$ then
\begin{eqnarray*}
&&\limsup_{\varepsilon\rightarrow0}\varepsilon\log
P(\sup_{0\le t\le1}|X_t^\varepsilon-X_n^\varepsilon(t)\ge\delta_0|)\\
&\le& \limsup_{\varepsilon\rightarrow0}\ \varepsilon\log(
P(\tau_n^\varepsilon\ge1,T_n^\varepsilon\le1)+P(\tau_n^\varepsilon\le1))\\
&\le& \limsup_{\varepsilon\rightarrow0}\ \varepsilon\log P(\tau_n^\varepsilon\ge1,T_n^\varepsilon\le 1)\vee
\limsup_{\varepsilon\rightarrow0}\varepsilon\log P(\tau_n^\varepsilon\le1)\\
&\le&(C-\int_0^{\delta_0^2}\frac{ds}{\eta(s)+\rho})\vee(-\frac{\delta}{4d\sup_{k\le n}A_{n,k}}).
\end{eqnarray*}

By letting $n\rightarrow\infty$, we have
$$
\lim_{n\rightarrow\infty}\limsup_{\varepsilon\rightarrow0}\ \varepsilon\log
P(\sup_{0\le
t\le1}|X_t^\varepsilon-X_n^\varepsilon(t)\ge\delta_0|)\le
C-\int_0^{\delta_0^2}\frac{ds}{\eta(s)+\rho}.
$$

Since $\rho\rightarrow 0$ as $\delta\rightarrow0,$ taking limit with
$\delta,$ the right hand side of the inequality tends to $-\infty $.
We complete the proof. $\square$

\section{Proof of Lemma \ref{l2}}\vskip0.2in

Define the test function
$$\varphi_{\rho}(x):=\exp{\big(\int_0^x\frac{ds}{\eta_R(s)+\rho}\big)}.$$

Let $$Y_n^l(t):=F_n(l)(t)-F(l)(t),\quad Z_n^l(t):=|Y_n^l(t)|^2.$$

For any $l$ with $e(l)\le\alpha$ and $\delta>0$ small enough (less than $1$), define

$$\tau_n(l):=\inf\{t\ge 0,|Y_n^l(t)|>\delta\}$$
and
$$\tau_R:=\inf\{t,|F_n(l)(t)|\vee|F(l)(t)|>R\}.$$

As in the proof of Lemma \ref{l1}, by (\ref{keji1}), it's clear that $\tau_R\rightarrow\infty$ as $R\rightarrow\infty$.

Denote $$e_s:=b(s,F_n(l)(\frac{[ns]}{n}))-b(s,F(l)(s))$$
$$h_s:=\sigma(s,F_n(l)(\frac{[ns]}{n}))-\sigma(s,F(l)(s))$$

Then by the chain rule, we have
$$
\aligned\varphi_{\rho}( Z_n^l(t\wedge\tau_n(l)\wedge\tau_R)
)&=1+2\int_0^{t\wedge\tau_n(l)\wedge\tau_R}\varphi'_{\rho}( Z_n^l(s) )
\langle Y_n^l(s),e_s\rangle ds\\
&\quad+2\int_0^{t\wedge\tau_n(l)\wedge\tau_R}\varphi'_{\rho}( Z_n^l(s) ) \langle
Y_n^l(s),h_sl'(s)\rangle
ds .\endaligned $$

Since $\sigma, b$ are bounded with $x$ for any fixed $t$, it follows that
$$\aligned |F_n(l)(t)-F_n(l)(\frac{[nt]}{n})|&\le\int_{\frac{[nt]}{n}}^t|b(s,F_n(l)(\frac{[ns]}{n}))|ds
+\int_{\frac{[nt]}{n}}^t||\sigma(s,F_n(l)(\frac{[ns]}{n}))|||l'(s)|ds\\&\le
C_\alpha\frac{1}{\sqrt{n}},\endaligned$$

\noindent where $0<C_\alpha$ is independent of $n.$ We have used H\"{o}lder inequality in the last step. As in proof of Lemma \ref{l1},
for $s\le t\wedge\tau_n(l)\wedge\tau_R$, we have
$$\aligned\langle
Y_n^l(s),e_s\rangle
&=\langle Y_n^l(s),b(s,F_n(l)(\frac{[ns]}{n}))-b(s,F_n(l)(s))\rangle\\&\quad
+\langle Y_n^l(s),b(s,F_n(l)(s)-b(s,F(l)(s))\rangle
\\&\le \delta G(s)H(\frac{C_\alpha}{\sqrt{n}})+f(s)\eta(Z_n^l(s))\\&
\le(f(s)+G(s))(\eta(Z_n^l(s))+\delta H(\frac{C_\alpha}{\sqrt{n}})).\endaligned $$

Similarly, we have
$$\aligned\langle
Y_n^l(s),h_s\rangle
\le (f(s)+G(s))(\eta(Z_n^l(s))+\delta H(\frac{C_\alpha}{\sqrt{n}})).\endaligned
$$

Take $\rho_n=\delta H(\frac{C_\alpha}{\sqrt{n}}).$ It follows that
$$
\aligned\varphi_{\rho_n}( Z_n^l(t\wedge\tau_n(l)\wedge\tau_R))\le
1+2\int_0^{t\wedge\tau_n(l)\wedge\tau_R} \varphi_{\rho_n}(Z_n^l(s))(f(s)+G(s))(1+|l'(s)|) ds .\endaligned $$

Since
$$\int_0^t(f(s)+G(s))|l'(s)|ds\le\Big(\int_0^t(f(s)+G(s))^2ds\Big)^\frac{1}{2}
\Big(\int_0^t|l'(s)|^2ds\Big)^\frac{1}{2}<\infty,$$
 by Gronwall's lemma and letting $R\rightarrow\infty,$ we have
$$\varphi_{\rho_n}( Z_n^l(1\wedge\tau_n(l)))\le\exp\{2\int_0^1
(f(s)+G(s))(1+|l'(s)|)ds\}. $$

Take supremum with $l\in\{l:e(l)\le\alpha\}$ and let $n\rightarrow\infty$, since $\varphi_\rho(x)$ is
increasing in $x,$ it follows that
$$
\limsup_{n\rightarrow\infty}\varphi_{\rho_n}(
\sup_{l:e(l)\le\alpha} Z_n^l(1\wedge\tau_n(l)))\le \exp\{2\int_0^1
(f(s)+G(s))ds+2\sqrt{2\alpha}\Big(\int_0^t(f(s)+G(s))^2ds\Big)^\frac{1}{2}\}.
$$

Since $\{\tau_n(l)>1\}=\{\sup_{0\le t\le1}|F_n(l)(t)-F(l)(t)|\le\delta\},$
we only need to show $\tau_n(l)>1$ for all $l\in\{l:e(l)\le\alpha\}$
and $n$ large enough. If not, there exists $\delta>0,$ a subsequence
$\{n_k,k\ge 1\}$ of positive integers and
$l_{n_k}\in\{l:e(l)\le\alpha\}$ such that $\tau_{n_k}(l_{n_k})\le
1.$ Then
$$\aligned
\varphi_{\rho_{n_k}}(\delta^2)&=\varphi_{\rho_{n_k}}(
Z_{n_k}^{l_{n_k}}(1\wedge\tau_{n_k}(l_{n_k})))\\&
\le\exp\{2\int_0^1
(f(s)+G(s))ds+2\sqrt{2\alpha}\Big(\int_0^t(f(s)+G(s))^2ds\Big)^\frac{1}{2}\}<\infty. \endaligned$$

Let $k\rightarrow\infty$. Then by the definition of $\varphi_\rho,$ the left hand side tends to $\infty$.
This is a contradiction. So we have
$$\sup_{l:e(l)\le\alpha}\sup_{0\le
t\le1}|F_n(l)(t)-F(l)(t)|\le\delta$$

\noindent for any $\delta>0$ for $n$ large enough. We complete the proof. $\square$

\noindent\textbf{Proof of Theorem \ref{Euler}}\quad Let $F_n$ and
$X_n^\varepsilon$ are Euler approximation of $F$ and $X^\varepsilon$
 respectively with the same scale. Notice that
 $X_n^\varepsilon(s)=F_n(\sqrt{\varepsilon}B)(s),$ where $B$ is the Brownian motion. It's
 clear that $F_n$ is continuous for each $n.$ Since it's well known that $\sqrt{\varepsilon}B$
 satisfies large deviation principle, now according to Lemma \ref{l1} and Lemma \ref{l2} and
 Theorem 4.2.23 in \cite{DZ}, we know that our process $X^\varepsilon$ also satisfies large
 deviation principle. We complete the proof. $\square$

\section{Large deviation principle in the general case}\vskip0.2in

\noindent\textbf{Proof of Lemma \ref{l3}}\quad Define
$\xi^\varepsilon(t):=|X^\varepsilon(t)|^2 $, and
$$
\varphi (x):=\exp{\big(\lambda\int_0^x\frac{ds}{\gamma(s)+1}\big)}.
$$

Define $\tau_R^\varepsilon:=\inf\{t>0,\xi^\varepsilon(t)\geq R^2\}$, it's clear that
$\tau_R\rightarrow\infty$ as $R\rightarrow\infty$ since the solution is non explosive 
under condition (\ref{zengzhang}). Then by It\^o's formula and taking expectation on 
both sides, it follows that
$$
\aligned\mathbb{E} \varphi(\xi^\varepsilon(t\wedge\tau_R^\varepsilon)) &\le 1+\mathbb{E}
\int_0^{t\wedge\tau_R^\varepsilon}\varphi'(\xi^\varepsilon(s))
g(s)(\gamma(\xi^\varepsilon(s))+1)ds\\&
\quad+2\varepsilon\mathbb{E}\int_0^{t\wedge\tau_R^\varepsilon}\varphi''(
\xi^\varepsilon(s)) g(s)(\gamma(\xi^\varepsilon(s))+1)^2ds
.\endaligned $$

By the definition of $\varphi$ and Gronwall's lemma, it follows that
$$ \mathbb{E}\varphi(\xi^\varepsilon(t\wedge\tau_R^\varepsilon))\le
e^{T\int_0^tg(s)ds}\varphi(|x_0|^2)$$

\noindent where $T=2\lambda^2\varepsilon+\lambda$. Let $t=1.$ Then
$$
P(\tau_R^\varepsilon\le 1)\varphi(R^2)\le
e^{T\int_0^1g(s)ds}\varphi(|x_0|^2).
$$

That is,
$$ \varepsilon\log P(\tau_R^\varepsilon\le 1)\le
T\varepsilon\int_0^1g(s)ds-\lambda\varepsilon\int_{|x_0|^2}^{R^2}\frac{dx}{\gamma(x)+1}.
$$

Take $\lambda=\frac{1}{\varepsilon}$, and let $\varepsilon\rightarrow0$ and $R\rightarrow\infty$ 
subsequently. We complete the proof. \quad $\square$

\noindent\textbf{Proof of Lemma \ref{l4}}\quad Define
$Z^l(t):=|F(l)(t)|^2 $, and
$$
\varphi(x):=\exp{\big(\int_0^x\frac{ds}{\gamma(s)+1}\big)}.
$$

It follows by (\ref{zengzhang}) that,
$$
\aligned\varphi(Z^l(t))\le\varphi(|x_0|^2)+\int_0^tg(s)(1+2|l'(s)|)\varphi(Z^l(s))ds.\endaligned $$

Now by Gronwall's lemma, it follows that
$$
\aligned\varphi(Z^l(t))&\le\varphi(|x_0|^2)\exp\{\int_0^tg(s)(1+2|l'(s)|)ds\}\\&
\le\varphi(|x_0|^2)\exp\{(1+2\sqrt{\alpha})\sqrt{\int_0^1g^2(s)ds}\}\endaligned $$
holds for any $l\in \{l,e(l)\le\alpha\}$ and $0\le t\le 1$. We have used H\"older inequality 
in the last inequality. Since $\varphi$ is increasing, it follows that
$$
\aligned\varphi(\sup_{e(l)\le\alpha}\sup_{0\le t\le 1}Z^l(t))
\le\varphi(|x_0|^2)\exp\{(1+2\sqrt{\alpha})\sqrt{\int_0^1g^2(s)ds}\}.\endaligned $$

By the definition of $\varphi$, it follows that
$\sup_{e(l)\le\alpha}\sup_{0\le t\le 1}Z^l(t)<\infty.\qquad\square$

In what follows, we will consider large deviation principle of
solutions of stochastic differential equation (\ref{SDE}) without the
integrable condition (\ref{keji1}) on $\sigma$ and $b$. To this end, we will
use the bounded approximation method.

For $R>0,$ define $$m_R(t)=\sup_{|x|\le R}\{|b(t,x)|\vee||\sigma(t,x)||\}.$$
Then $m_R\in L^2([0,t])$  for any $t>0$ by (\ref{keji}). Let 
$$b_i^R(t,x):=(-m_R(t)-1)\vee b_i(t,x)\wedge(m_R(t)+1),$$
$$\sigma_{ij}^R(t,x):=(-m_R(t)-1)\vee \sigma_{ij}(t,x)\wedge(m_R(t)+1)$$ and
$$b_R(t,x):=(b_1^R(t,x),\cdots,b_d^R(t,x)),\quad\sigma_R(t,x):=(\sigma_{ij}^R)(t,x).$$

Then
$$b_R(t,x)=b(t,x),\ \sigma_R(t,x)=\sigma(t,x),\quad \forall |x|\le R,\ t\in [0,1].$$

It's obvious that $b_R, \sigma_R$ satisfy (\ref{keji1}), and satisfy (\ref{lianxu1}) with
the same $f$ and $\eta_R$. Let $X_R^\varepsilon $ be the solution of
$$
dX^\varepsilon_R(t)=\varepsilon^{\frac{1}{2}}\sigma_R(t,X^\varepsilon_R(t))dB_t+b_R(t,X^\varepsilon_R(t))dt,\
X_R^\varepsilon(0)=x_0.$$

For $l\in C_{0}([0,1],R^m)$ with $e(l)<\infty $, let $F_R(l)(t) $ be the solution of
\beqn\label{bijin}
dF_R(l)(t)=\sigma_R(t,F_R(l)(t))l'(t)dt+b_R(t,F_R(l)(t))dt,\
F_R(l)(0)=x_0.
\eeqn

If $\sup_{0\le t\le 1}|F(l)(t)|\le R ,$ where $F(l)$ is the solution of differential equation
(\ref{fuzhu}), then $F(l)$ is the solution of differential equation
(\ref{bijin}). By uniqueness of solutions, one can see that $F_R(l)(t)=F(l)(t)$ for $0\le t\le 1$.
Define
$$I_R(f)=\inf\{\frac{1}{2}e(l):F_R(l)=f\},\ f\in C_{x_0}([0,1],R^d).
$$

Then
\[I_R(f)= I(f),f\in C_{x_0}([0,1],R^d),\quad\sup_{0\le t\le1}|f(t)|\le R.\]

By Theorem \ref{Euler}, $\{\mu_\varepsilon^R,\varepsilon>0\}$ satisfies a large deviation principle (with $I_R(\cdot)$).

\noindent\textbf{Proof of Lemma \ref{l5}}\quad By Lemma \ref{l4},
for $\alpha>0$, there exists $R>0$ such that
\beqn\label{youjie1}\sup_{l:e(l)\le\alpha}\sup_{0\le t\le1}|F(l)(t)|\le R.\eeqn

Thus
$$F_R(l)=F(l),\quad \forall l\in\{l:e(l)\le\alpha\}.$$

Let $\{f_n\}$ be a sequence in $\{f:I(f)\le\alpha\}.$ Then there
exists a sequence $l_n\in C_{0}([0,1],R^m)$ such that $F(l_n)=f_n$ and
$\frac{1}{2}e(l_n)\le \alpha+\frac{1}{n}.$ So there exists a limit point
$l\in C_{0}([0,1],R^m)$ of $\{l_n\}$ such that $\frac{1}{2}e(l)\le \alpha$. According to
(\ref{youjie1}) we have $F_R(l_n)=F(l_n)=f_n$, and $F_R(l_n)$
converges uniformly (over $[0,1]$) to $F_R(l)=F(l)$ (up to a subsequence). Let
$f=\lim\limits_{n\rightarrow\infty}f_n=F(l)$. Then
$I(f)=\frac{1}{2}e(l)\le \alpha$. So $\{f:I(f)\le\alpha\}$ is
compact. $\square$

\noindent\textbf{Proof of Theorem \ref{wujie}}\quad Repeat the proof of theorem $E$ in Fang and Zhang \cite{FZ} word by word,
we can prove that Theorem \ref{wujie} holds under conditions (\ref{lianxu1}) and (\ref{zengzhang}), so we omit it here.\quad $\square$

\end{document}